%% file: root.tex
\documentclass[letterpaper, 10 pt, conference]{ieeeconf} 

\IEEEoverridecommandlockouts
\usepackage{amsmath,amsfonts,amssymb,amsthm}
\usepackage{graphicx}
\usepackage{xcolor}
\usepackage{comment}
\usepackage{algorithm}
\usepackage{algpseudocode}
\usepackage{mathtools}
\usepackage{subcaption}
\usepackage{dsfont}
\usepackage[font=footnotesize,labelfont=bf]{caption}
\usepackage{lettrine}
\usepackage{cite}

\newcommand{\dsqbL}{[\kern-0.15em[}
\newcommand{\dsqbR}{]\kern-0.15em]}

\input{lpsymbols}

\input{mymacros} 
\newcommand{\bvarepsilon}{\boldsymbol{\varepsilon}}

\graphicspath{ {./fig/} }

\title{\bf Planning Strategies for Lane Reversals in Transportation Networks
}
\author{Salom\'{o}n Wollenstein-Betech, Ioannis Ch. Paschalidis, and Christos G. Cassandras 
\thanks{This work was supported in part by NSF under grants ECCS-1509084, DMS-1664644, CNS-1645681, IIS-1914792, and CMMI-1454737, by AFOSR under grant FA9550-19-1-0158, by ARPA-E under grant DE-AR0001282, by the MathWorks, by the ONR under grant N00014-19-1-2571, and by the NIH under grant 1R01GM135930.}
\thanks{The authors are with the
Division of Systems Engineering, Boston University, Brookline, MA, 02446, USA {\tt\small \{salomonw, cgc, yannisp\}@bu.edu}}
}

\begin{document}
    \maketitle
    \input{sections/00_Abstract.tex}

\input{sections/01_Intro.tex}
    \input{sections/02_Problem-Formulation.tex}
    \input{sections/03_Lane-Assignment}

\input{sections/04_LASO-opt}
    \input{sections/05_Numerical-Results}

    \input{sections/06_Conclusion}
    \small
    \bibliographystyle{IEEEtran}
    \bibliography{ref}
\end{document}

%% file: lpsymbols.tex
\def\bd{{\mathbf d}}

\def\bt{{\mathbf t}}

\def\bw{{\mathbf w}}
\def\bx{{\mathbf x}}  

\def\bz{{\mathbf z}}

\def\bN{{\mathbf N}}

\def\texitem#1{\par\smallskip\noindent\hangindent 25pt
               \hbox to 25pt {\hss #1 ~}\ignorespaces}

\newcommand{\scrA}{\mathcal{A}}

\newcommand{\scrF}{\mathcal{F}}
\newcommand{\scrG}{\mathcal{G}}

\newcommand{\scrV}{\mathcal{V}}
\newcommand{\scrW}{\mathcal{W}}

\newcommand{\bbeta}{\boldsymbol{\beta}}

\newcommand{\btheta}{\boldsymbol{\theta}}

\newcommand{\bpsi}{\boldsymbol{\psi}}


%% file: mymacros.tex
\newcommand{\be}[1]{\begin{equation}\label{#1}}
\newcommand{\benon}{\begin{equation*}}  
\newcommand{\bemuln}[1]{\begin{multline}\label{#1}}
\newcommand{\bemul}{\begin{multline*}}
\newcommand{\bea}{\begin{eqnarray*}}
\newcommand{\eea}{\end{eqnarray*}}
\newcommand{\been}[1]{\begin{eqnarray}\label{#1}}
\newcommand{\eeen}{\end{eqnarray}}
\newcommand{\began}[1]{\begin{gather}\label{#1}}
\newcommand{\bega}{\begin{gather*}}
\newcommand{\bealn}[1]{\begin{align}\label{#1}}
\newcommand{\beal}{\begin{align*}}
\newcommand{\bealatn}[2]{\begin{alignat}{#1}\label{#2}}
\newcommand{\bealat}{\begin{alignat*}}
\newcommand{\bexalatn}[1]{\begin{xalignat}\label{#1}}
\newcommand{\bexalat}{\begin{xalignat*}}

\newcommand{\mbb}{\mathbb}

\theoremstyle{plain} \newtheorem{thm}{Theorem}[section]

\newtheorem{rem}{Remark}



%% file: sections/00_Abstract.tex
\begin{abstract}

This paper studies strategies to optimize the lane configuration of a transportation network for a given set of Origin-Destination demands using a planning macroscopic network flow model. 
The lane reversal problem is, in general, NP-hard since the optimization is made over integer variables. 
To overcome this burden, we reformulate the problem using a piecewise affine approximation of the travel latency function which allows us to exploit the total unimodularity property of Integer Linear Programming (ILP). 
Consequently, we transform the ILP problem to a linear program by relaxing the integer variables. 
In addition, our method is capable of solving the problem for a desired number of lane reversals which serves to perform cost-benefit analysis. We perform a case study using the transportation network of Eastern Massachusetts (EMA) and we test our method against the original lane configuration and a projected lower bound solution. 
Our empirical results quantify the travel time savings for different levels of demand intensity. We observe reduction in travel times up to 40\% for certain links in the network.

\end{abstract}


%% file: sections/01_Intro.tex
\section{Introduction} \label{sec:intro}

\lettrine{T}{raffic} congestion and urban mobility are among the most relevant topics regarding the development and sustainability of cities. 
Thanks to the increasing ability to gather, communicate, and process data, we are closer towards the establishment of ``Smart Cities''.
Thus, urban-policy decisions today will have a considerable impact on how urban and suburban mobility unfolds in the future. 
Therefore, it is imperative to think of short, medium, and long term solutions to reduce traffic congestion.

An indisputable solution to ease traffic congestion is to increase the network's capacity. One way to accomplish this without the economical and societal cost of building new roads, is to dynamically reverse the direction of certain lanes on specific roads of the network. This strategy is known as \textit{Contraflow Lane Reversals} (or simply \emph{Lane Reversals})~\cite{zhou1993intelligent,xue2000intelligent,hausknecht2011dynamic} and has already been implemented during rush hour times in many cities (e.g., Mexico City, Montreal) in order to mitigate traffic congestion.
\begin{figure}[t]
    \centering
    \begin{subfigure}{0.55\columnwidth}
        \includegraphics[trim = 10em 8em 12em 0em ,clip,width=1\columnwidth]{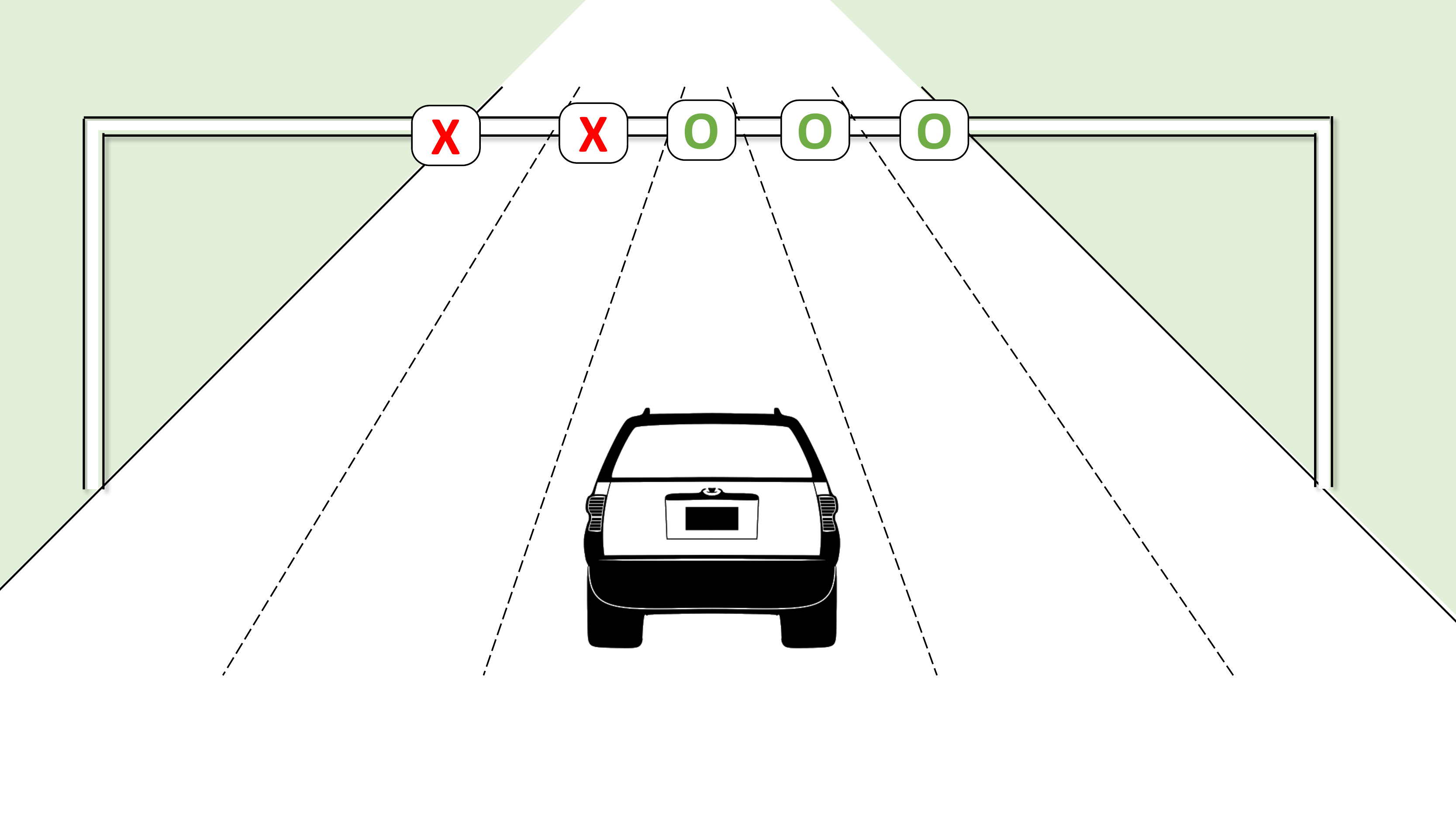}
        \vspace{-1em}
        \caption{}
        \label{fig:contraflow-signal}
    \end{subfigure}    
    \begin{subfigure}{0.42\columnwidth}
        \includegraphics[width=1\columnwidth]{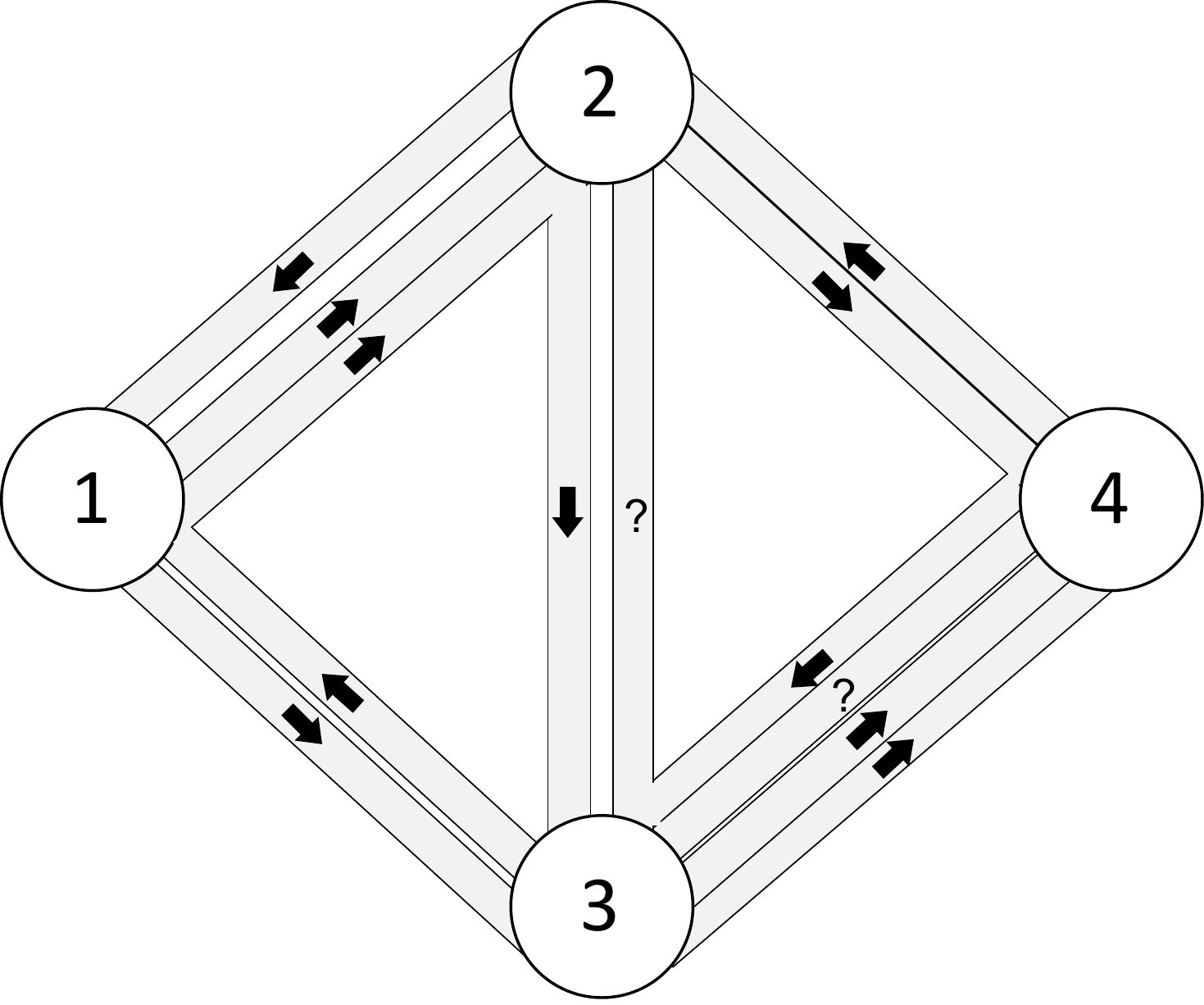}
        \caption{}
        \label{fig:contraflow-diagram}
    \end{subfigure}
    \caption{(a) Illustration of a typical lane reversal signal. (b) Diagram of the Braess' network with specific number of lanes and directions.}
    \vspace{-1.5em}
\end{figure}

Still, current contraflow lane reversals are used in practice in limited situations, particularly when emergency evacuations are needed (e.g., Hurricane Florence in South Carolina~\cite{fausset2018hurricane}), or after the culmination of a big event (e.g., NASCAR races in New Hampshire~\cite{nh2019contraflow}). 
In both cases, the reversal must be carefully planned before being implemented. 

Presumably, the main challenge when applying lane reversals is to effectively communicate the updated direction of the lanes to drivers.
Current systems utilize overhead signals indicating the direction of the lane (see Fig.~\ref{fig:contraflow-signal}), or they simply define structured schedules throughout the day. 
However, only a few cities have adopted this technology, and even then, it has only been done for a small fraction of the lanes.
In part this is because at the time when the reversal occurs, a gridlock can be produced by poor human ability to respond to such a coordinated change.
In addition, lane reversals may produce more congestion due to increases in traffic accidents~\cite{martinez2021contraflowaccident}. 
Thus, the spread of contraflow lane reversals has been curtailed because of the lack of effective means for communicating the status of the transportation infrastructure to drivers.

Fortunately, the rapid development of Connected and Automated Vehicles (CAVs) can address such limitations. 
The promising ability of CAVs to communicate with the infrastructure offers the possibility to implement contraflow lane reversals more aggressively either by doing it for more roads, or by dynamically changing the direction of a single road more regularly.  

In this paper we take a macroscopic planning view to address the lane reversal problem. Rather than investigating when to reverse a single lane, we focus on 
selecting the best lane configuration that alleviates the  traffic congestion when the overall demand pattern of travelers is  considered. Figure~\ref{fig:contraflow-diagram} shows a diagram illustrating the problem. 

Most of the research solving the contraflow lane reversal problem has focused either on reversing lanes for evacuation routing plans during emergencies, or to ease traffic congestion. 
For evacuation route planning, the problem has been solved using simulation and network flow models. 
Simulation methods, \cite{jha2004emergency} and~\cite{theodoulou2004alternative}, showed that evacuation route capacity can be improved by 53\% and 73\%, respectively, when designing appropriate lane reversals. 
Furthermore, for network flow models,~\cite{cova2003network} proposes a mixed-integer programming formulation for which a solution is found using a generic solver. 
They report improvements on total population evacuation time of the order of 30\% to 40\% for the Salt Lake City network.
In addition to these numerical results~\cite{kim2008contraflow} showed that this network flow problem is NP-hard and provided a greedy heuristic algorithm. 

The research concerned with reducing traffic congestion (and not evacuation planning) has focused in reversing lanes in single bottleneck roads (typically tunnels or bridges) or selecting the reversing lanes considering the full transportation network.
For single roads, rule-based~\cite{ampountolas2019motorway} and fuzzy~\cite{zhou1993intelligent,xue2000intelligent} controllers have been proposed and they typically rely on the fundamental diagram of traffic flow~\cite{lighthill1955kinematic}.
For network-wide settings, its canonical mathematical representation is an Integer Linear Program (ILP) which is NP-hard.
To tackle it,~\cite{chu2019dynamic} uses a distributed alternating direction method of multipliers (ADMM) to decompose the problem into smaller integer programs. Their numerical results report improvements in travel times of 61\% in New York City. 
In addition,~\cite{hausknecht2011dynamic,meng2008optimizing} use genetic algorithms and report increases in efficiencies of 72\% for the city of Austin.
Alternatively,~\cite{levin2016cell} solves a microscopic cell-transmission version of the problem using a heuristic based on congestion estimates. Their results show a 21\% reduction in total system travel time.
Note that none of the models above provide guarantees to find an optimal solution. However, their numerical results show that these algorithms work well in practice and motivate the use of lane reversals to alleviate traffic congestion. 

The contribution of this paper is threefold. First, and our main result, is to transform the NP-hard ILP lane reversal problem with fixed flows to a tractable, polynomial-time, linear program (LP). To that end, we approximate the convex function of the lane reversal problem with a piecewise affine and convex function. Interestingly, this approximation maintains the total unimodularity of the ILP. Therefore, we can solve the problem using any linear programming method and it is ensured to recover an integer solution. Second, we propose an algorithm that calculates a lower-bound of the flow-aware problem using the original convex function. Finally, we test our methods and provide numerical results using the Eastern Massachusetts (EMA) transportation network. 

The remainder of the paper is organized as follows: 
In Section~\ref{sec:model} we present the  preliminaries and problem formulation.
In Section~\ref{sec:lane-assigment} we introduce the piecewise affine approximation to reduce the NP-hard convex integer programming (IP) problem to a tractable LP. 
In Section~\ref{sec:LASO} we discuss the optimality conditions of the lane reversal planning problem and propose an algorithm to obtain a lower bound.
In Section~\ref{sec:Numerical-results} we show the performance of the proposed method over a case study using the EMA network. Finally, we draw some conclusions in Section~\ref{sec:conclusion}. 
In this paper, all vectors are column vectors and denoted by bold lowercase letters. We use ``prime'' to denote transpose, and use $\mathbf{1}$ to denote the vector whose entries are all ones.

%% file: sections/02_Problem-Formulation.tex
\section{Model and Problem Formulation} \label{sec:model}

Consider a transportation network as a strongly-connected directed graph $\scrG( \scrV, \scrA )$, where $\scrV$ is the set of nodes (intersections) and $\scrA$ is the set of arcs (roads).
For every arc $(i,j)$ in $\scrA$, we denote the number of lanes assigned to it with $z_{ij}$ and its total capacity (veh/h) with $m_{ij}=z_{ij}c_{ij}$ where $c_{ij}$ is the capacity per lane of arc $(i,j)$. Moreover, we let $n_{ij}=n_{ji}$ be the sum of the number of lanes in arcs $(i,j)$ and $(j,i)$.
Therefore the maximum number of lanes that can be assigned to $(i,j)$ or $(j,i)$ is $n_{ij}$.
Let the node-arc incidence matrix of $\scrG$ be $\bN \in {\{ {0,1, - 1} \}^{| \scrV | \times | \scrA|}}$ and $\textbf{e}_{ij}\in \mbb{R}^{|\scrA|}$ be a zero vector with an entry equal to $1$ corresponding to arc $(i,j)$. 
Let $\bw_k = (w_{sk},w_{tk})$ denote an Origin-Destination (OD) pair, the total number of OD pairs be $K$, and 
$\scrW = \{{\bw_k}:{\bw_k} = ( {{w_{sk}},{w_{tk}}} ), \,k =1,\dots,K \}$ be the set of all OD pairs.  
For every $k=1,\dots,K$, let ${d^{\bw_k}} \ge 0$ be the demand of flow (veh/hr) that travels from origin $w_{sk}$ to destination $w_{tk}$, and denote by $\bd^{\bw_k} \in \mathbb{R}^{|\scrV|}$ the vector of all zeros except for the coordinates of nodes $w_{sk}$ and $w_{tk}$ which take values $-d^{\bw_k}$ and $d^{\bw_k}$, respectively. 
Let $x_{ij}$ be the total flow on arc $(i,j)$ and $\bx=(x_{ij};\ i,j\in\scrA)$ be the vector of all arc flows in the transportation network. 

To maintain demand satisfaction and conservation of flow, we let an arc flow vector $\bx$ be feasible if $\bx\in\scrF$ where 
\begin{align} \label{eq:feasibility}
    \scrF = \Big\{  \bx \in \mathbb{R}_+^{|\scrA|} \ | \ \bx =
		\sum\limits_{k=1}^K {{\bx^{\bw_k}}},\,  
		& \bN{\bx^{\bw_k}}  = {\bd^{\bw_k}}, \hspace{2em} \\[-1em] 
	    & \hspace{1.5em} k=1,\dots,K  \Big\}, \notag
\end{align}
and $\bx^{\bw_k}$ is the flow vector attributed to OD pair $\bw_k$.

In addition, let $t_{ij}(x_{ij}, z_{ij}) : \left\{\mathbb{R}_{\geq0}, \mathbb{N}_{\geq0}\right\} \mapsto \mathbb{R}_{\geq0}$ be the \emph{latency cost} (i.e., travel time) function for arc $(i,j)$ which depends on the arc's flow and on the number of lanes assigned to the arc. Using the same structure as in~\cite{Beckmann1955}, we characterize $t_{ij}(x_{ij}, z_{ij})$ as:
\begin{align} \label{eq:latency function}
t_{ij}(x_{ij}, z_{ij}) = t_{ij}^0f\Big(\frac{x_{ij}}{c_{ij}z_{ij}}\Big),
\end{align}
where $f(\cdot)$ is a strictly increasing, positive, and continuously differentiable function, and $t_{ij}^0$ is the free-flow travel time on arc $(i,j)$. 
We set $f(0)=1$, which ensures that if there is no
constraint on the arc's capacity, the travel time $t_{ij}$ equals the free-flow travel time $t^0_{ij}$ . 
To explicitly characterize $f(\cdot)$ we employ the widely-used \emph{bureau of public roads} (BPR) function~\cite{BPR}:
\begin{equation} \label{eq:bpr}
f(x)=1+0.15x^4.
\end{equation}
However, we could use any other polynomial function as discussed in~\cite{WollensteinSunEtAl2019}. Finally, we let the vector of travel latency functions be $\bt(\bx,\bz)=(t_{ij}(x_{ij}, z_{ij}); \forall i,j\in\scrV)$.

Using these definitions, we now can formulate the \emph{Lane-Assignment System-Optimal Traffic Assignment Problem} (LASO-TAP) as follows:
\begin{subequations} \label{eq:LASO-TAP}
\begin{align}
    \min_{\bx,\bz} \quad & 
    \bt(\bx,\bz)'\bx \label{eq:LASO-TAP-obj} \\
    \text{s.t.} \quad & z_{ij} + z_{ji} \leq n_{ij}, \quad \forall (i,j)\in\scrA, \label{eq:LASO-TAP-cnstr-max-lanes} \\
    & \bx\in\scrF, \quad \bz\in\mathbb{N}^{|\scrA|}_{+}, \label{eq:LASO-TAP-cnstr-feasible set}
\end{align}
\end{subequations}
where in the objective~\eqref{eq:LASO-TAP-obj} we are minimizing the overall travel time, in constraint~\eqref{eq:LASO-TAP-cnstr-max-lanes} we ensure that the number of assigned lanes does not exceed the maximum number of available lanes, and in constraint~\eqref{eq:LASO-TAP-cnstr-feasible set} we restrict $\bx$ and $\bz$ to be feasible vectors which complies with demand satisfaction and conservation of flow (c.f.~\eqref{eq:feasibility}). 

The LASO-TAP is difficult to solve for several reasons. First, the interaction of the decision variables in~\eqref{eq:LASO-TAP-obj} makes the objective non-convex. This is because when multiplying~\eqref{eq:latency function} by $x_{ij}$ we get  $\gamma_{ij} x_{ij}^5/z_{ij}^4$ where $\gamma_{ij}=0.15 t^0_{ij}/c^4_{ij}$.
Second, we are optimizing over a set of integer variables $\bz$ which makes the problem computationally intractable for large instances of the problem. In the following, we present modifications to this problem and propose efficient approximate solutions for it. 

\begin{rem}
    {\normalfont 
    Note that the LASO-TAP as stated in~\eqref{eq:LASO-TAP} is seeking a System-Optimal (SO) solution to the Traffic Assignment Problem (TAP). We can expect a SO behavior when vehicles collaborate with each other when deciding their routes (which is investigated for the presence of CAVs in cities~\cite{houshmand2019penetration,WollensteinHoushmandEtAl2020}). However, by slightly modifying the objective function, the same algorithms that solve a SO TAP are capable of solving the User-Centric (UC) TAP~\cite{Patriksson1994a}. Hence, in this manuscript we focus on the SO case, but our framework can easily accommodate both SO and UC TAPs.
    }
\end{rem}

%% file: sections/03_Lane-Assignment.tex
\section{Lane Assignment Problem} \label{sec:lane-assigment}


Recall that in the literature (e.g.,~\cite{kim2008contraflow}) the LASO problem has been identified as NP-hard. Thus, let us consider exclusively the \emph{Lane Assignment} (LA) problem by fixing the arc flows $\bx$. This is to
\begin{subequations} \label{eq:LA}
\begin{align}
    \min_{\bz} \quad & 
    J(\bz):=\bt(\bx,\bz)'\bx \label{eq:LA-obj} \\
    \text{s.t.} \quad & z_{ij} + z_{ji} \leq n_{ij}, \quad \forall (i,j)\in\scrA, \label{eq:LA-max-lanes} \\
    & \bz\in\mathbb{N}^{|\scrA|}_{+}, \label{eq:LA-cnstr-feasible set}
\end{align}
\end{subequations}
where we have fixed $\bx$ and therefore $J(\bz)$ in~\eqref{eq:LA-obj} only depends on $\bz$.
The problem LA, as stated in~\eqref{eq:LA}, is an integer program (IP) with a convex objective and linear constraints. 
The convexity comes from the fact that each element $J_{ij}:=t(x_{ij}, z_{ij})x_{ij}$ is convex in $z_{ij}$ when we use the BPR function defined in~\eqref{eq:bpr}. 
This is because its second derivative is nonnegative for all $z_{ij}\geq0$, and by the fact that the sum of convex functions is convex. Besides the objective, a key property of this formulation is that the left-hand side of the constraint matrix formed by~\eqref{eq:LA-max-lanes} is \emph{unimodular} (all submatrices have determinant $1$ or $-1$), and the right-hand side is an integer-valued vector $\mathbf{n}$.

Even with these structural properties, solving~\eqref{eq:LA} is computationally expensive as it requires optimizing over integer variables. Now, we present our main result which exploits convexity, monotonicity and total unimodularity to provide an optimal solution to~\eqref{eq:LA} using LP. 

\subsection{Piecewise affine approximation} \label{subsec:LA-piecewise-affine}
\begin{figure}
    \centering
    \includegraphics[width=.95\columnwidth]{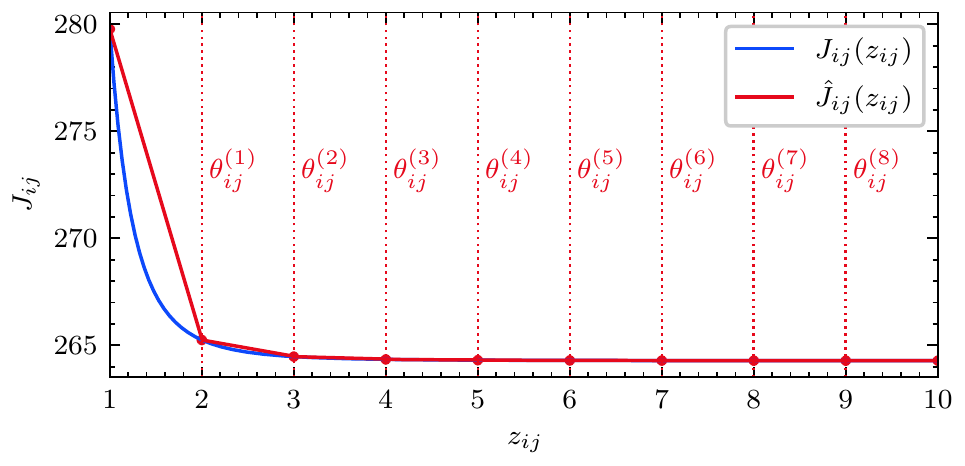}
    \caption{Piecewise affine approximation of the cost function for a specific arc in the network.}
    \label{fig:approx}
    \vspace{-1.7em}
\end{figure}

In order to obtain an optimal solution to the LA problem, we approximate the function $J_{ij}$ with a piecewise affine and convex function, as shown in Fig.~\ref{fig:approx}. 
Let $\btheta_{ij}=(1,\dots,n_{ij})$ be a vector of \emph{thresholds} for every arc $(i,j)$ and let its $k$th entry be $\theta_{ij}^{(k)}$. These thresholds serve as breakpoints of the piecewise affine approximation and we have intentionally selected them to be integers as we are interested in evaluating the function at these points. For each segment $(\theta^{(k)}_{ij}, \theta^{(k+1)}_{ij})$ we let $\beta^{(k)}_{ij}$ be the slope of the corresponding affine segment, and we will use $\beta^0_{ij}$ to refer to the function's intercept. 

For every arc $(i,j)$ we let $\bvarepsilon_{ij}\in\{0,1\}^{n_{ij}}$ be a vector of $n_{ij}$ \emph{slack variables}, each corresponding to an assignment of a lane to arc $(i,j)$. Hence, we have that $z_{ij}=\sum_{k=1}^{n_{ij}}\varepsilon^{(k)}_{ij}$. With these definitions, we rewrite the LA problem as
\begin{subequations} \label{eq:approx}
\begin{align}
    \min_{\bvarepsilon} \quad & \hat{J}(\bvarepsilon) := \sum_{(i,j)\in\scrA} \beta^0_{ij} + \bbeta_{ij}' \bvarepsilon_{ij}  \label{eq:approx-obj} \\
    \text{s.t.} \quad & \sum_{k=1}^{n_{ij}} \big( \varepsilon^{(k)}_{ij} + \varepsilon^{(k)}_{ji} \big) \leq n_{ij}, \quad \forall (i,j)\in\scrA, \label{eq:approx-cnstr-max-lanes} \\
    & \sum_{k=1}^{n_{ij}} \varepsilon^{(k)}_{ij} \leq  \theta^{(l)}_{ij}, \quad \hspace{-0.5em} \forall (i,j)\in\scrA, \ l=1,\dots,n_{ij},  \label{eq:approx-cnst-epsilon}\\
    & \bvarepsilon_{ij}\in\{0,1\}^{n_{ij}}, \quad \forall (i,j)\in\scrA, \label{eq:approx-cnst-nonneg}
\end{align}
\end{subequations}
where $\bvarepsilon=(\bvarepsilon_{ij}; \ \forall (i,j)\in\scrA)$. 
Using the model in~\eqref{eq:approx}, we state the main results of our paper.

\begin{thm}[Equivalence of Problems~\eqref{eq:LA} and~\eqref{eq:approx}]
    If the intersections of the piecewise affine segments defining $\hat{J}_{ij}$ coincide with the value of the original function $J_{ij}$ and $J_{ij}$ is convex and monotonically decreasing, then, solving problem~\eqref{eq:approx} is equivalent to solving problem~\eqref{eq:LA}.
\end{thm}
\begin{proof}
    We prove it by construction. Observe that the main difference between~\eqref{eq:LA} and~\eqref{eq:approx} is the use of an approximated piecewise affine objective function $\hat{J}$ rather than $J$. 
    We argue that since~\eqref{eq:LA} is over the integer-valued variables $z_{ij}$, the only points for which $J_{ij}$ is evaluated are the ones which coincide with the breakpoints $\btheta_{ij}$ (see Figure~\ref{fig:approx}) and are exactly evaluated by our assumption on $\hat{J}_{ij}$.
    Moreover, by the monotonicity and convexity of $J$ (which implies that $\beta^{(k)}_{ij}\leq \beta^{(k+1)}_{ij}$) we have that if $z_{ij}$ is integer, then $\epsilon^{(k)}_{ij}=1$ if and only if all $\epsilon^{(l)}_{ij}=1$ for $l=1,\dots,k-1$. This allows the objective to be $\bbeta_{ij}'\bvarepsilon_{ij}$.
    Note that with this modification both objective functions are the same, and the constraints~\eqref{eq:LA-max-lanes} and~\eqref{eq:approx-cnstr-max-lanes} are identical by using the fact that $z_{ij}=\sum_{k=1}^{n_{ij}}\varepsilon^{(k)}_{ij}$. 
    Finally, constraint~\eqref{eq:approx-cnst-epsilon} restricts $\varepsilon^{(k)}_{ij}$ to be either $0$ or $1$. 
\end{proof}

\begin{thm}[Total unimodularity]
    If the integer variables $\bvarepsilon$ in~\eqref{eq:approx} are relaxed to be continuous, then, the solution to the relaxed version of~\eqref{eq:approx} is an integer solution.
\end{thm}
\begin{proof}
    First, observe that problem~\eqref{eq:approx} seeks to minimize a linear objective over a set of constraints whose matrix is unimodular (integer matrix). Next, notice that the right hand side vector of~\eqref{eq:approx-cnst-nonneg} and~\eqref{eq:approx-cnst-epsilon} is integer-valued since $n_{ij}$ and $\theta^{(k)}_{ij}$ are integers for all $(i,j)$ and all $k=1,\dots,n_{ij}$. Thus, by using the Unimodularity Theorem~\cite[Thm. 11.15]{AhujaMagnantiEtAl1993}, we can now relax the integer variables $\bvarepsilon$ to be continuous, and the resulting LP is certified to recover an integer solution. 
\end{proof}

Interestingly, these two results imply that the LA problem can be solved with any polynomial-time LP algorithm. 
The only disadvantage of~\eqref{eq:approx} with respect to~\eqref{eq:LA} is that it requires $\sum_{(i,j)\in\scrA} n_{ij}$ decision variables instead of $|\scrA|$, i.e., the total number of lanes instead of the number of arcs.
Nonetheless, because these extra variables grow linearly with respect to the original ones, the computational cost is not too high. 
In addition, note that solving the LA problem obtains an upper bound to the optimal value of the LASO-TAP (since LA selects the best $\bz$ but LASO-TAP optimizes both $\bx$ and $\bz$).

The problem as stated in~\eqref{eq:approx} is a lane \emph{assignment} problem rather than a lane \emph{reversal} problem. 
This distinction is negligible when we are capable to reverse any number of lanes in the network. 
In reality, the transportation infrastructure is not very flexible and might not be able to handle many lane reversals.
Hence, we are interested in considering a \emph{sparse} lane assignment problem in which we limit the number of lane reversals.
To achieve this, let $\bz^0\in\mathbb{N}^{|\scrA|}_{\geq0}$ be a preferred lane configuration, for example, the current network condition. 
Then, our formulation~\eqref{eq:approx} can handle the sparse case by adding an extra term in the objective function that penalizes an $\ell_1-$norm deviation from the nominal $\bz^0$. i.e., for every arc, we would penalize $\|\mathbf{1}'\bvarepsilon_{ij}-\bz^0_{ij}\|_1$. 
It is straightforward to see that these problems can be reformulated as LPs~\cite{bertsimas1997introduction}.
Conversely, we can also include this term as a constraint (rather than a cost) and fix a maximum number of lane reversals allowed in the optimization procedure.

\subsection{Lower Bound} \label{subsec:LA-lower-bound}

To compare our solution with a theoretical benchmark, we design a lower bound solution to the problem. 
A natural approach to generate lower bounds when dealing with IPs is to relax the integer variables to continuous ones to obtain a convex programming problem. 
We follow this idea and derive optimality conditions for the relaxed problem. To do so, we quantify the impact on the objective when we reverse a small fraction of a lane while the flows $\bx$ remain unchanged. 
Formally, we perturb an infinitesimal amount of capacity $\delta$, and estimate its impact on $J$ as follows:
\begin{align}
    \frac{\partial J}{\partial z_{ij}} = & \lim_{\delta\xrightarrow[]{}0} \Big(\big(x_{ij}t_{ij}(x_{ij},z_{ij}+\delta)+x_{ji}t_{ji}(x_{ji},z_{ji}-\delta)\big) \notag\\
    & \quad - \big(x_{ij}t_{ij}(x_{ij},z_{ij})+x_{ji}t_{ji}(x_{ji},z_{ji})\big)\Big)/\delta. \label{eq:partial-derivatives} 
\end{align}

Note that the way in which we are calculating these derivatives assumes that the flow vector $\bx$ remains constant. This is aligned with solving the LA instead of the LASO problem and we expect it to be a good estimate even for the LASO when $\delta$ is sufficiently small. 
Then, an optimality condition of the relaxed problem is that for all $(i,j)$ in $\scrA$,
\begin{align*}
    &\partial J/\partial z_{ij}=0, \quad \text{for } z_{ij}\in(0,n_{ij}), \\
    &\partial J/\partial z_{ij}\geq0, \quad \text{for } z_{ij}=0, \\
    &\partial J/\partial z_{ij}\leq0,  \quad \text{for } z_{ij}=n_{ij}.
\end{align*}
That is, there is no benefit to reversing $\delta$ units of capacity for any arc $(i,j)$. 
Thus, since the relaxed problem is convex, we propose a \emph{feasible direction} method described in Alg.~\ref{alg:feasible-direction} to solve the $\xi$-optimal relaxed LA problem where we use $\Pi$ to indicate a projection to the closest feasible value. 
\begin{algorithm}[t] 
\caption{Relaxed LA problem} \label{alg:feasible-direction}
\begin{algorithmic}[1]
    \While {$\| \nabla_{\bz}J \| \geq \xi$}
        \For {$(i,j) \in \scrA$}
            \State $\frac{\partial J}{\partial z_{ij}} \leftarrow \text{using}~\eqref{eq:partial-derivatives}$
            \State $z_{ij} \leftarrow \Pi(z_{ij} - \alpha_{ij} \frac{\partial J}{\partial z_{ij}})$
            \State $\alpha_{ij}\leftarrow$ update coordinate step-size
        \EndFor
    \EndWhile
\end{algorithmic}
\end{algorithm}

%% file: sections/04_LASO-opt.tex
\section{LASO Optimization} \label{sec:LASO}

So far, we have discussed a simpler version of our LASO-TAP formulation in~\eqref{eq:LASO-TAP} by fixing the vector flow $\bx$ in~\eqref{eq:LA}. Now, let us discuss the optimality conditions of the LASO-TAP (i.e., when we also optimize the flows) as well as some algorithms to find good solutions for it. 

\subsection{Optimality conditions} \label{subsec:LASO-optimality-conditions}

Similar to~\eqref{eq:partial-derivatives}, the idea is to measure the impact on the objective when we reverse a single lane in the network. 
Consider reversing a lane in arc $(i,j)$; then its impact on the objective can be calculated with:
\begin{subequations}
\begin{align}
    \psi_{ij}=&\big(x_{ij}^{+1}t_{ij}(x_{ij}^{+1},z_{ij}+1)+x_{ji}^{+1}t_{ji}(x_{ji}^{+1},z_{ji}-1)\big) \notag \\
    & - \big(x_{ij}t_{ij}(x_{ij},z_{ij})+x_{ji}t_{ji}(x_{ji},z_{ji})\big) \notag \\
    & + \sum_{\mathclap{(k,l)\in\scrA/\{(i,j), (j,i)\}}}\Big( \big(x^{+1}_{kl}t_{kl}(x_{kl}^{+1},z_{kl})+x_{lk}^{+1}t_{lk}(x_{lk}^{+1},z_{lk})\big) \notag \\
    & - \big(x_{kl}t_{kl}(x_{kl},z_{kl})+x_{lk}t_{lk}(x_{lk},z_{lk}\big)\Big), \notag 
\end{align}
\end{subequations}
where $x^{+1}_{ij}$ indicates the updated flows when we reverse a lane on arc $(i,j)$. 
Then, the necessary optimality conditions are $\psi_{ij}\geq0$ for $z_{ij}=0,\dots,n_{ij}-1$ and $\psi_{ij}\leq0$ for $z_{ij}=n_{ij}$, for all $(i,j)$ in $\scrA$.
Note that to calculate this, we have to solve a TAP to get $x^{+1}_{ij}$.
Therefore, to evaluate the full vector $\bpsi$, we are required to solve $\left|\scrA\right|$ TAPs. 
Checking this condition in a sequential algorithm is computationally expensive, but could be used regularly as a measure of closeness to optimality.

\subsection{Algorithms} \label{subsec:LASO-algorithms}
Since the conditions presented above are expensive to evaluate, we propose a faster approximation by assuming that routing decisions do not change when a single lane is reversed, i.e., $\bx$ remains fixed. In this case, the \emph{relaxed optimality condition} is the same as in Section~\ref{subsec:LASO-optimality-conditions} but we use an approximated $\hat{\bpsi}$ instead of $\bpsi$, defined as follows:
\begin{subequations}
\begin{align}
    \hat{\psi}_{ij} = &\big(x_{ij}t_{ij}(x_{ij},z_{ij}+1)+x_{ji}t_{ji}(x_{ji},z_{ji}-1)\big) \notag \\
    & \qquad - \big(x_{ij}t_{ij}(x_{ij},z_{ij})+x_{ji}t_{ji}(x_{ji},z_{ji})\big). \label{eq:derivative-estimates}
\end{align}
\end{subequations}
In other words, reversing the direction of a lane does not improve the cumulative travel times for that road. 
Note that it is hard to establish that $\hat{\psi}_{ij}$ approximates $\psi_{ij}$ in expectation as we anticipate routing decisions to fluctuate when we change a lane's direction. This is especially true for congested lanes. 
However, the idea is that by computing $\hat{\psi}_{ij}$ fast, we can design algorithms that yield good solutions for the LASO-TAP.

We list three algorithms that can be employed to solve the LASO-TAP: $(i)$ Use the classical \emph{Frank-Wolfe} algorithm used to solve the TAP~\cite{Patriksson1994a} and at each step include an update on the capacities using the derivative estimates~\eqref{eq:derivative-estimates}. In this manner, the optimization will be trying to minimize~\eqref{eq:LASO-TAP-obj} by updating both $\bx$ and $\bz$ at every step.
$(ii)$ Utilize a \emph{sequential} approach consisting of solving the TAP, fixing $\bx$, and then taking a descent step in the negative direction of~\eqref{eq:partial-derivatives}. 
$(iii)$ Employ a \emph{fully sequential} approach which repeatedly solves the TAP and then solves the LA problem~\eqref{eq:approx}. Due to space limitations of this manuscript, we leave the development of these approaches for future work. 

%% file: sections/05_Numerical-Results.tex
\section{Numerical Results and Case Studies} \label{sec:Numerical-results}
\begin{figure}
    \centering
    \includegraphics[trim={3cm 3.4cm 3cm 3.5cm},clip, width=0.65\columnwidth]{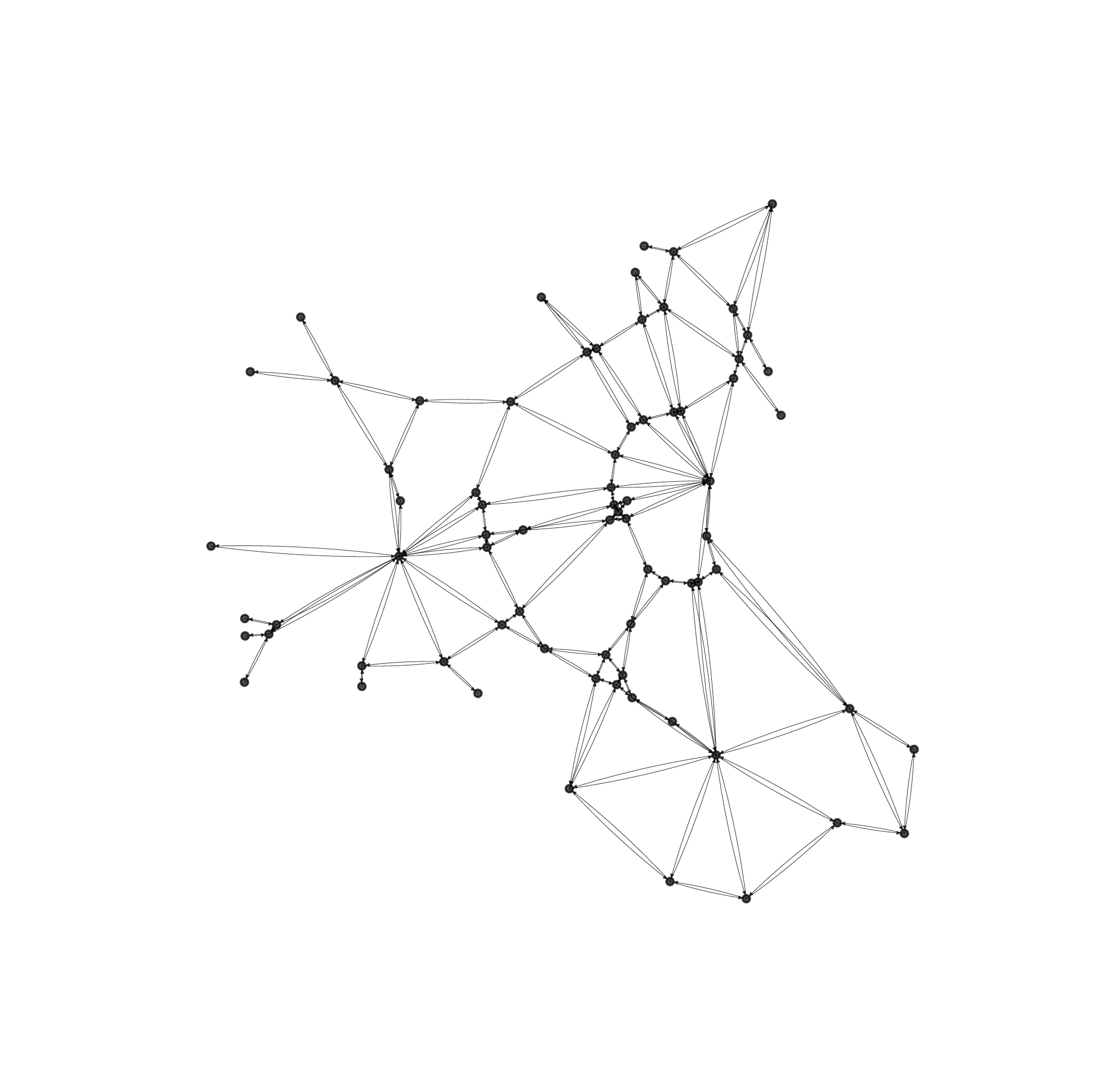}
    \caption{EMA transportation network composed of 74 nodes , 258 arcs, 581 lanes, and 1113 OD pairs.}
    \label{fig:EMA-net}
    \vspace{-1.5em}
\end{figure}
To test the effectiveness of our method we designed two case studies using the EMA transportation network. The EMA network (Fig. \ref{fig:EMA-net}) consists of 74 nodes, 258 links, 581 lanes, and 1113 OD pairs, and it captures mobility in the context of suburban/urban topologies where we expect lane reversal strategies to be beneficial.
This is because arterial roads typically have a large number of lanes and, at the same time, they experience high traffic congestion.
The values of the network topology (e.g., capacities, free flow speeds), and the code used to perform the experiments are publicly available in our online repository.\footnote{https://github.com/salomonw/contraflow-lane-reversal}

\subsection{Performance of the LA solution}

\begin{figure}
    \centering
    \includegraphics[width=\columnwidth]{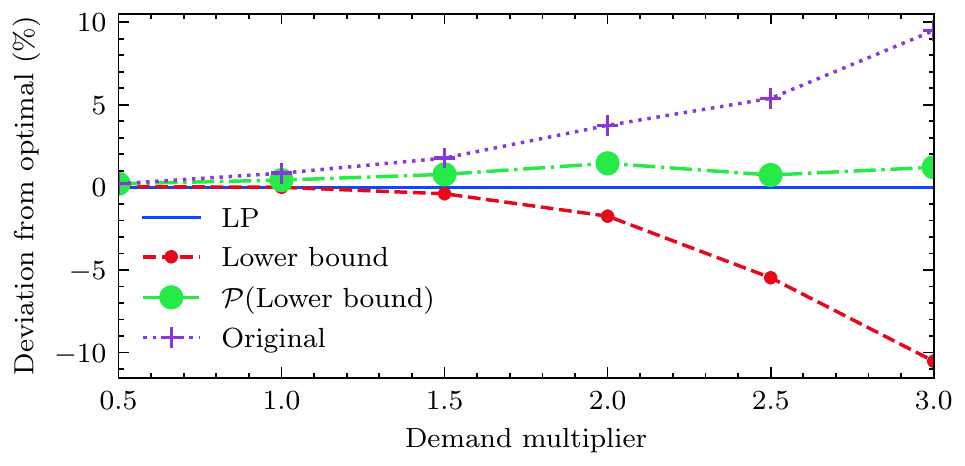}
    \caption{Deviation in percentage of the optimal travel times achieved when solving the LA problem~\eqref{eq:approx} against the \emph{Original} lane configuration and the projected relaxation solution, \emph{$\mathcal{P}$(Lower Bound)}, for different congestion levels. 
    For example, for a demand multiplier equal to 2.5, the overall travel time for the original configuration is around 5\% higher than when lanes are rearranged using~\eqref{eq:approx}. 
    In the same spirit, we show a lower bound obtained when solving Alg.~\ref{alg:feasible-direction}, and an upper-bound when projecting this solution to its closest integer.
    As expected, we see that the higher the congestion, the higher the benefit of utilizing lane reversals.}
    \label{fig:peformance}
    \vspace{-0.5em}
\end{figure}

\begin{figure}
    \centering
    \includegraphics[width=\columnwidth]{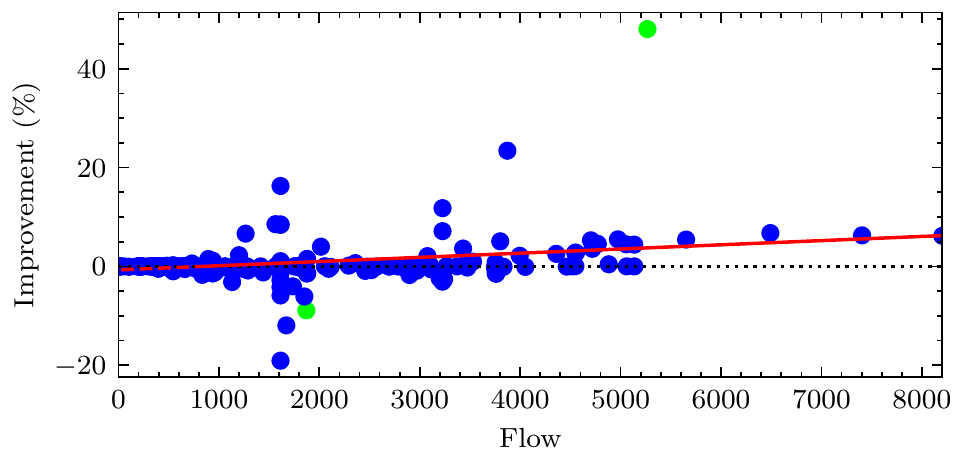}
    \caption{Each dot represents an arc and reports its flow and its travel time improvement with respect to the nominal capacity before reversing the best 30 lanes. We observe travel time improvements of up to 40\%. In addition, we see that, on average, the higher the flow in a link, the higher the benefit we can achieve. Note that we have arcs that experience higher travel times (those below the dotted line). This happens because as we improve travel times when adding capacity to link $(i,j)$, we automatically decrease the capacity of $(j,i)$ (the green dots show this for a particular pair of arcs with opposite directions).
    }
    \label{fig:link-performance}
    \vspace{-1.6em}
\end{figure}

We are interested in comparing the performance of the solution of~\eqref{eq:approx} against other approaches and for different congestion levels. 
In particular, we compare its solution with $(i)$ the \emph{original} lane configuration of the transportation network and $(ii)$ projecting the solution of the relaxed lower bound presented in Section~\ref{subsec:LA-lower-bound} into integers. 
That is, rounding each continuous variable obtained when running Alg.~\ref{alg:feasible-direction} to its closest integer. 
Figure~\ref{fig:peformance} shows the performance of these approaches for different demand levels. 
The results show that for low congestion levels, the achievable benefits are marginal. 
That is because the network is not congested and therefore there are no benefits to rearranging the lane directions. 
As we begin to increase demand, we see that we can achieve overall travel time improvements of almost 10\% when compared with the original lane configuration, and around 1\% when compared to the projected relaxation approach.
It is important to stress that these results are measuring the performance of the overall travel time of commuters.
Hence, when we analyze the performance of each individual arc, we see improvements of up to 40\%. Figure~\ref{fig:link-performance} shows the performance per arc relative to its flow.
In addition to the comparison with these two benchmarks, we also include the lower bound value achieved by following Alg.~\ref{alg:feasible-direction} and its projected version.
It is interesting to see that the projected lower bound solution is close to the optimal LA one. 
This is encourages further research to develop projected solutions for the optimal relaxed LASO-TAP.

\subsection{Maximum number of reversals}

Our next experiment explores the Pareto optimal frontier using the \emph{sparse} version of the LA problem described in the last paragraph of Section~\ref{sec:lane-assigment}. 
That is, we explore the performance as we increase the maximum number of allowable reversals. 
To achieve this, we include a constraint that fixes the maximum number of lane reversals and we solve the problem for different values.
Figure~\ref{fig:max-reversals} shows this relationship when the demand multiplier is equal to $1.5$.
As expected, the first lanes contribute the most to the improvement of the overall travel times. 
These results are of particular interest to urban planners that  need to prioritize the most critical roads. 
Moreover, our results imply that it is not necessary to apply too many lane reversals to achieve a large fraction of possible improvements. 
For example, by making between 15 to 20 reversals (out of 581 possible lane reversals) in Fig.~\ref{fig:max-reversals}, the solution has already reached a good performance level. 


%% file: sections/06_Conclusion.tex
\section{Conclusion}\label{sec:conclusion}
Most existing approaches to tackle the contraflow lane reversal problem use standard heuristic algorithms to solve the Integer Programming (IP) problem encountered. 
In this paper, instead of using these techniques, we approximate the objective using a piecewise affine function that reformulates the problem as a Linear Program (LP). We show that this LP obtains the same optimal solution as the IP problem when the flows are fixed. 
In addition to this technical result, our empirical results show the potential benefits of using contraflow lane reversals for a case study in the Eastern Massachusetts Area. Our numerical results show that lane reversals can be used as a tool to reduce commuter travel times and that they provide good performance even when the total number of lane reversals is restricted.

\begin{figure}
    \centering
    \includegraphics[width=\columnwidth]{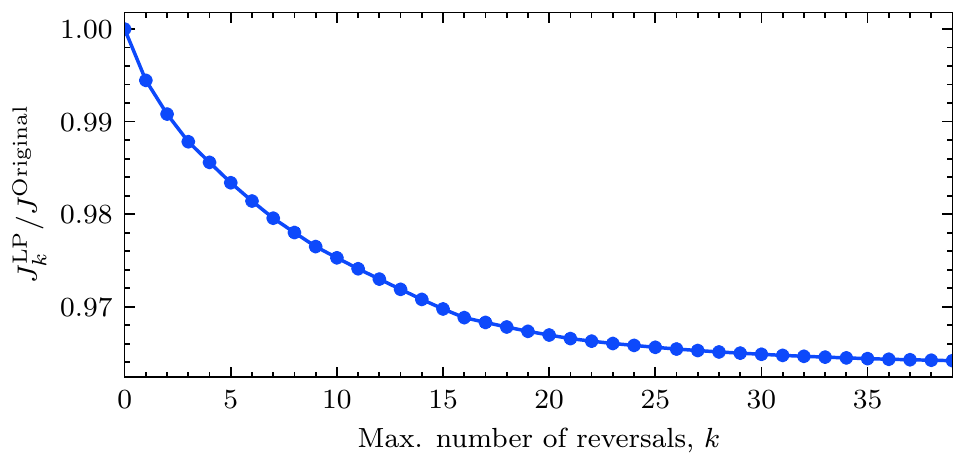}
    \caption{Performance on the overall commuting travel times as the number of contraflow lane reversals increases.}
    \label{fig:max-reversals}
    \vspace{-1.8em}
\end{figure}